\documentclass[11pt]{amsart}
\usepackage{amsthm, amsmath}
\usepackage{amssymb} 
\usepackage{amscd}

\theoremstyle{plain}
\newtheorem{thm}[subsection]{Theorem}
\newtheorem{lem}[subsection]{Lemma}
\newtheorem{prop}[subsection]{Proposition}

\theoremstyle{remark}
\newtheorem{rem}[subsection]{Remark}
\def\s{\subsection{\kern-5pt}}

\def\rond{\kern 1pt{\scriptstyle\circ}\kern 1pt}

\def\pr{\noindent\emph{Proof} : }

\def\gal{\mathrm{Gal}(\bar k/k)}
\def\lr#1{\langle {#1}\rangle}
\def\bb{\Bbb}

\def\R{\bb R}
\def\C{\bb C}
\def\P{\bb P}

\def\F{\bb F}
\def\Sym{\mathrm{Sym}}

\def\qfl#1{\buildrel {#1}\over {\longrightarrow}}

\def\iso{\vbox{\hbox to .8cm{\hfill{$\scriptstyle\sim$}\hfill}
\nointerlineskip\hbox to .8cm{{\hfill$\longrightarrow $\hfill}} }}

\begin{document}
\title[Jacobians among abelian threefolds: a geometric approach]{Jacobians among abelian threefolds:\\ a geometric approach}
\author[Arnaud Beauville]{Arnaud Beauville}
\address{Laboratoire J.-A. Dieudonn\'e\\
UMR 6621 du CNRS\\
Universit\'e de Nice\\
Parc Valrose\\
F-06108 Nice cedex 2, France}
\email{arnaud.beauville@unice.fr}
 \author[Christophe Ritzenthaler]{Christophe Ritzenthaler}
 \address{Institut de Math\'ematiques de Luminy\\
Universit\'e Aix-Marseille - CNRS\\
Luminy Case 907, F-13288 Marseille cedex 9, France}
 \email{ritzenth@iml.univ-mrs.fr}
\date{\today}
 
\begin{abstract}
Let $(A,\theta )$ be a principally polarized abelian threefold over a perfect field $k$, not isomorphic to a product over $\bar k$. There exists a canonical extension $k'/k$, of degree $\leq 2$, such that $(A,\theta )$ becomes isomorphic  to a Jacobian over $k'$. The aim of this note is to give a geometric construction of this extension.
\end{abstract}

\maketitle

\section{Introduction}

\par Let $(A,\theta )$ be a principally polarized abelian variety of dimension 3 over a field $k$. If $k$ is algebraically closed, $(A,\theta )$ is the Jacobian variety of a curve $C$ (or a product of Jacobians). If $k$ is an arbitrary perfect field the situation is more subtle (see Proposition \ref{serre} below): there is still a curve $C$ defined over $k$, but either $(A,\theta )$ is 
 isomorphic to $JC$, or  they become isomorphic only after a quadratic  extension $k'$ of $k$, uniquely determined by  $(A,\theta )$. 

\par Now given  $(A,\theta )$, how can we decide if it is a Jacobian, and more precisely determine the extension $k'/k$? For $k\subset \C$, a solution is given in \cite{LRZ} in terms of modular forms. Here we propose a geometric approach, based on a construction of Recillas. We have to assume that $A$ admits a rational theta divisor $\Theta $, and a rational point $a\in A(k)$ outside a certain explicit divisor $\Sigma_A \subset A$. This guarantees that the curve $\Theta \cap (\Theta +a )$ is smooth, and that there exists $b\in A(k)$ such that the involution $z\mapsto b-z$ acts freely on that curve. The quotient $X_a $ is a non-hyperelliptic genus 4 curve; its canonical model lies on a unique quadric $Q\subset \P^3$. Then for $\mathrm{char}(k)\neq 2$ the extension $k'$ is $k(\sqrt{\mathrm{disc}(Q)})$ (we will give more detailed statements in \S 3).
\par  The proof has two steps. We consider first the case where $(A,\theta )$ is a Jacobian, and prove that in that case the quadric $Q$ is  $k$-isomorphic to $\P^1\times \P^1$ (\S 2). Then we treat the case 
where $(A,\theta )$ is \emph{not} a Jacobian, and prove that the nontrivial automorphism of the extension $k'/k$ exchanges the two rulings of $Q$ (\S 3); this is enough to prove the theorem.
\bigskip
\section{Recillas' construction}
\par  Throughout the paper we work over a perfect field $k$. \par  In this section we fix a non-hyperelliptic curve $C$ of genus $3$ (that is, a smooth plane quartic curve), defined over $k$. We will denote by $K$ its canonical class.
We assume that the principal polarization of $JC$ can be defined by a theta divisor $\Theta $ \emph{defined over $k$}. There exists a degree 2 divisor class $D$ such that $\Theta $ is the image of 
 $\Sym^2C-D$ in $JC$; this class is unique, hence defined over $k$. 
 Note that since $C$ is not hyperelliptic, $\Theta $ is smooth and the map $E\mapsto E-D$ induces an isomorphism of $\Sym^2C$ onto $\Theta $. 
 \par  We choose a point $a\in JC(k)$; we are interested in the curve  $\tilde X_a :=\allowbreak\Theta \cap (\Theta +a )$.
  Put $b=K+a-2D\in JC(k)$; we have  $-\Theta =\Theta +a-b$. The involution $z\mapsto b-z$ exchanges $\Theta $ and $\Theta +a$, hence induces an involution $\iota $ of  $\tilde X_a$. We define a divisor $\Sigma _{JC}\subset JC$ by $\Sigma_{JC} =\Sigma'_{JC} \cup \Sigma''_{JC}$, where
  \[\Sigma' _{JC}=\{2E-K\ |\ E\in \Sym^2C\}\qquad\hbox{and}\qquad \Sigma_{JC}'' =C-C\ .
 \]
  
\begin{prop} The curve $\tilde X_a :=\Theta \cap (\Theta +a )$ is smooth and connected if and only if $a \in JC\smallsetminus \Sigma_{JC} $. If this is the case, the involution $\iota $ of $\tilde X_a$ is fixed point free. 
\label{rec}
\end{prop}
\pr  Throughout the paper it  will be convenient to use the following nota\-tion: given a divisor class $d$ of degree 2 on $C$ with $h^0(\mathcal{O}_C(d))=1$, we  denote by $\lr{d}\in\Sym^2C$ the unique effective divisor in the class $d$.

\par   Let $z\in \tilde X_a$. By \cite{K}, thm. 2, the tangent space $\P T_z(\Theta )\subset$ $ \P T_z(JC)=\P^2$ is identified with the line spanned by the divisor $\lr{D+z}\in \Sym^2C$. Similarly $\P T_z(\Theta +a )$ is identified with the line spanned by $\lr{D+z-a} \in \Sym^2C$;  the intersection $\tilde X_a$ is singular at $z$ if and only if these two lines coincide. If this happens,  then either
\par  $\bullet$ the two divisors $\lr{D+z}$ and $\lr{D+z-a} $ have a common point, which implies $a \in C- C$; or
\par  $\bullet$ $\lr{D+z}+\lr{D+z-a}\sim K$, which implies $a \in \Sigma'_{JC}$.
\par  Conversely, if $a \in\Sigma'_{JC} $, we have $K+a\sim 2E$ with $E\in \Sym^2C$; then
$z=E-D$ is a singular point of  $\tilde X_a$. If $a\sim p-q$, with $p,q\in C$, the intersection $ \tilde X_a$ is reducible, equal to $(C+p-D)\cup (K-D-q-C)$. 
 
\par  Assume now $a \notin \Sigma_{JC} $, so that $\tilde X_a$ is smooth;  the exact sequence
\begin{eqnarray*}
0\rightarrow \mathcal{O}_{JC}(-\Theta -(\Theta +a))\longrightarrow \mathcal{O}_{JC}(-\Theta)\oplus \mathcal{O}_{JC}(-(\Theta+a)) \longrightarrow\hfill \\
\hfill\longrightarrow \mathcal{O}_{JC}\longrightarrow \mathcal{O}_{\tilde X_a}\rightarrow 0
\end{eqnarray*}
gives $h^0(\mathcal{O}_{\tilde X_a})=1$, hence $\tilde X_a$ is connected.
If a  point  $z\in \tilde X_a$ is fixed by $\iota $ it satisfies $2(D+z)\sim K+a$, which implies $a\in\Sigma'_{JC} $.\qed \medskip

\s We assume from now on $a\notin\Sigma_{JC} $. We denote by $X_a$ the quotient curve $\tilde X_a/\iota $. The adjunction formula gives
\[K_{\tilde X_a}\sim (\Theta +(\Theta +a))^{}_{|\tilde X_a}=2\Theta ^3=12\ ,\quad\hbox{hence}\quad  g(\tilde X_a)=7\ \hbox{ and }\ g(X_a)=4\ .
\]
If $\mathrm{char}(k)\neq 2$, the principally polarized abelian variety $JC$ is canonically isomorphic to the Prym variety associated to 
the double covering $\tilde X_a\rightarrow X_a$ (\cite{B}, \S 3.b; \cite{R}). 

\par  We embed $\tilde X_a$ into $\Sym^2C\times \Sym^2C$ by $z\mapsto (\lr{D+z}, \lr{K+a-(D+z)})$. Then   $\tilde X_a$ is identified with $s^{-1}(|K+a|)$, where $s:\Sym^2C\times \Sym^2C\rightarrow \Sym^4C$ is the sum map. The involution $\iota $ is induced by the involution of $\Sym^2C\times\allowbreak \Sym^2C$ which exchanges the factors. The map $s:\tilde X_a\rightarrow |K+a|$ factors through $\iota $, hence induces a 3-to-1 map $t:X_a\rightarrow |K+a|\ (\cong \P^1)$.  
The fibre of $t$ above $E\in|K+a|$ parametrizes the decompositions $E=d+d'$, with $d,d'$ in $\Sym^2C$.

\par  We now consider the involution $(d,d')\mapsto (\lr{K-d},\lr{K-d'})$ of  \break$\Sym^2C\times \Sym^2C$; it maps  $\tilde X_a$ onto     $s^{-1}(|K-a|)=\tilde X_{-a}$ and  commutes with $\iota$, hence induces an isomorphism $X_a\iso X_{-a}$. By composition with the map $X_{-a}\rightarrow |K-a|$ defined above we obtain another degree 3 map $t':X_{a}\rightarrow |K-a|$.
\par  The maps $t$ and $t'$ are defined over $k$; they define two $g^1_3$ on $X_a$, that is, two
linear series of degree $3$ and projective dimension $1$, defined over $k$.
\begin{lem}\label{lem}
The two  $g^1_3$ defined by $t$ and $t'$ on $X_a$ are distinct. 
\end{lem}
\pr Let us first observe that the degree 4 morphism $f:C\rightarrow \P^1$ defined by the linear system $|K+a|$ is separable. If this is not the case, we have $\mathrm{char}(k)=2$ and $f$ factors as $C\qfl{F} C_1\qfl{g} \P^1$, where $C_1/k$ is the pull back of $C/k$ by the automorphism $\lambda \mapsto \lambda ^2$ of $k$, $F$ is the Frobenius $k$-morphism and $g$ is separable of degree 2 (see \cite{H}, IV.2). But then $C_1$ is hyperelliptic, hence also $C$.

\par   Assume that the two linear series are the same. By the previous observation there exists a  divisor $E=p+q+r+s$   in 
$|K+a|$ consisting of 4 distinct points.  There must exist  $E'\in |K-a|$ such that $t^{-1}(E')=t'^{-1}(E)$. This means that for each decomposition $E=d+d'$ with $d,d'$ in $\Sym^2C$, we have $E'=\lr{K-d}+\lr{K-d'}$. 
\par  Let us write $\lr{K-p-q}=p'+q'$ and $\lr{K-r-s}=r'+s'$, so that $E'=p'+q'+r'+s'$. 
We must have $E'=\lr{K-p-r}+\lr{K-q-s}$, so we can suppose $\lr{K-p-r}=p'+r'$. Then $K-p-p'\sim q+q'\sim r+r'$, which implies $r'=q$, $q'=r$.  But then we get $K-p-q-r\sim p'$ and $a\sim s-p'$, which contradicts the hypothesis $a\notin \Sigma_{JC} $.\qed
\medskip
\par  We can now conclude:
\begin{prop}\label{jac}
For $a\notin \Sigma_{JC} $, the genus $4$ curve $X_a$ is not hyperelliptic; the unique quadric $Q\subset \P^3$ containing its canonical model is smooth and split over $k$ (that is, isomorphic to $\P^1\times \P^1$ over $k$). 
\end{prop}
\pr Since $X_a$ admits a base point free $g^1_3$ it cannot be hyperelliptic (\cite{ACGH}, p.\ 13). Let us denote the two distinct $g^1_3$ of $X_a$ by $|E|$ and $|E'|$. We have $E+E'\sim K_{X_a}$; by the base-point free pencil trick,  the multiplication map $H^0(X_a,E)\otimes H^0(X_a,E')\rightarrow H^0(X_a,K_{X_a})$ is an isomorphism.   
Thus the canonical map of $X_a$ is the composition of $(t,t'):X_a\rightarrow \P^1\times \P^1$ and of the Segre embedding $\P^1\times \P^1\hookrightarrow \P^3$, so the unique quadric containing the image is isomorphic to $\P^1\times \P^1$.\qed
\begin{rem}
If $C$ is hyperelliptic, Proposition \ref{rec} still holds, with  essentially the same proof. However Lemma \ref{lem} fails: in fact, we have $t'=\sigma \rond t$, where $\sigma :|K+a|\iso |K-a|$ is induced by the hyperelliptic involution. Actually in that case $X_a$ has a unique $g^1_3$, at least  if $\mathrm{char}(k)\neq 2$. Indeed  $\Theta $ has a singular point, given by the $g^1_2$ of $C$; on the other hand $JC$ is isomorphic to the Prym variety of $\tilde X_a/X_a$.  
By \cite{M2}, \S 7, Thm. (c),  this happens if and only if $X_a$ admits a unique $g^1_3$.
\end{rem}

\s
 The divisor $\Sigma'_{JC}$  is equal to $\mathbf{2}_*\Delta $, where $\mathbf{2}$ is the endomorphism $z\mapsto 2z$ of $JC$ and $\Delta $ is any symmetric theta divisor; thus it
  can be defined on any absolutely indecomposable principally polarized abelian threefold $(A,\theta )$, with no reference to the isomorphism $A\iso JC$. 
 The same holds for $\Sigma ''_{JC}$ provided $\mathrm{char}(k)\neq 2$. Recall indeed
  that there is  a canonical linear system on $A$, denoted $|2\theta |$, which contains the double of each
symmetric theta divisor. Then:
\label{sigma}
\begin{lem}
If $\mathrm{char}(k)\neq 2$, the divisor $\Sigma ''_{JC}=C-C$ is the unique divisor in $|2\theta |$ with multiplicity $\geq 4$ at $0$. 
\end{lem}\vskip -8pt
 This is quite classical if $k=\C$, see \cite{GG}. We do not know whether it still holds when $\mathrm{char}(k)= 2$.

\pr  The difference map $C\times C\rightarrow C-C$ is an isomorphism outside the diagonal $\Delta $, and contracts $\Delta $ to $0$; therefore the multiplicity of $C-C$ at $0$ is $-\Delta ^2=4$.
\par Let us prove the unicity; we may assume $k=\bar k$. We denote by $|2\theta |_0$ the subspace of elements of $|2\theta |$ containing $0$. The multiplicity at $0$ of an element of $|2\theta |$ is even: this follows from the ``inverse formula" of \cite{M1}, p. 331. Thus we have a projective linear map $\tau :|2\theta |_0\dasharrow |\mathcal{O}_{\P^2}(2)|$ which associates 
to a divisor  its quadratic tangent cone at $0$. Since $\dim |2\theta |_0=6$ and $\dim |\mathcal{O}_{\P^2}(2)|=5$, it suffices to prove that $\tau $ is surjective. 
For each  $E\in \Sym^2C$, the divisor $(\Sym^2C-E)+(\Sym^2C-(K-E))$ belongs to $|2\theta |_0$; by \cite{K}, thm. 2,  its tangent cone at $0$ is twice the line in $\P^2$ spanned by $E$. Since the double lines span the space of conics, $\tau $ is surjective.\qed

\bigskip
\section{The main result}
\par In this section we fix a  principally polarized abelian threefold $(A,\theta )$ over $k$. We assume that it is  absolutely indecomposable, that is, $(A,\theta )$ is not isomorphic over $\bar k$ to a product of two principally polarized abelian varieties. It is equivalent to say  that the theta divisor of $A$ is irreducible (over $\bar k$), or that $(A,\theta )_{\bar k}$ is isomorphic to the Jacobian of a curve \cite{OU}. This does not imply that $(A,\theta )$ itself is a Jacobian; indeed we have \cite{S}:
\begin{prop}\label{serre}
There exists a curve $C$ over $k$ and a character \break $\varepsilon^{}_A : \gal\rightarrow \{\pm 1\}$, uniquely determined, such that $(A,\theta )$ is $k$-isomorphic to $JC$ twisted by $\varepsilon^{}_A $. If $C$ is hyperelliptic, $\varepsilon^{} _A$ is trivial. 
\end{prop}
\s  In more down-to-earth terms this means the following. Let $k'$ be the extension of $k$ defined by the character $\varepsilon^{} _A$. Then: 
\par  $\bullet\ $ if $\varepsilon^{} _A=1$ (that is, $k'=k$), $JC$ is  $k$-isomorphic to $(A,\theta )$. This is the case if $C$ is hyperelliptic.
\par  $\bullet\ $ if $\varepsilon^{} _A\neq 1$ (that is, $k'$ is a quadratic extension of $k)$,
$JC$ is  isomorphic to $(A,\theta )$ over $k'$ \emph{but not over} $k$. More precisely, let $\sigma $ be the nontrivial automorphism of $k'/k$; there exists an isomorphism $\varphi :(A,\theta )\rightarrow JC$ such that ${}^\sigma \varphi =-\varphi $.
\label{explain}
\s Our aim is to describe geometrically the character $\varepsilon^{}_A $. We will compare it to the character associated to a smooth quadric $Q\subset \P^3_k$ in the following way: such a quadric admits two rulings defined over $\bar k$, so the action of $\gal$ on these rulings provides a character $\varepsilon^{} _Q:\gal\rightarrow \{\pm 1\}$. We will describe this character in more concrete terms below.

\par  We define  the divisor $\Sigma_A=
\Sigma'_A\cup \Sigma'' _A $ on $A$ as in  \ref{sigma}: we put  $\Sigma' _A=\mathbf{2}_*\Delta $ for any symmetric theta divisor $\Delta $; if  $\mathrm{char}(k)\neq 2$,  $\Sigma''_A$ is the unique divisor in $|2\theta |$ with multiplicity $\geq 4$ at $0$. An alternative definition, which works in all characteristics, is as follows: we choose an isomorphism $A\iso JC$ over $k'$ and  put $\Sigma ''_{A}=\varphi ^{-1}(C-C)$.  Since $C-C$ is symmetric this definition does not depend on the choice of $\varphi $.
\par   We assume that $A$ admits a theta divisor $\Theta $ defined over $k$. If $\Theta $ is singular, $C$ is hyperelliptic, hence $A\cong JC$ by Proposition \ref{serre}. Thus we may assume that $\Theta $ is smooth.   
\par  We also assume that there exists  a  point $a\in A(k)$ outside $\Sigma_A $. The divisor $-\Theta $ is in the class of the polarization $\theta $, hence
there is a unique $b\in A(k)$ such that $(-\Theta )+b=\Theta +a$; the involution $z\mapsto b-z$ exchanges $\Theta $ and $\Theta +a$.
\begin{thm}
Let $X_a$ be the quotient of the curve $\Theta \cap (\Theta +a)$ by 
the involution $z\mapsto b-z$. Then $X_a$ is a smooth curve of genus $4$, non hyperelliptic. Its canonical model lies in a smooth quadric $Q\subset \P^3$, and we have $\varepsilon^{}_A =\varepsilon^{} _Q$. 
\end{thm}
\pr Following \ref{explain} we choose an isomorphism $\varphi :(A,\theta )\rightarrow JC$ defined over $k'$. It induces an isomorphism of $\tilde X_a$ onto the corresponding curve $\tilde X_{\varphi (a)}\subset JC$, hence of $X_a$ onto $X_{\varphi (a)}$.  By remark \ref{sigma} $\varphi $ maps $\Sigma_A $ onto  $\Sigma _{JC}$, thus $\varphi (a)\notin \Sigma _{JC}$; then Proposition \ref{jac} tells us that $X_a$ is not hyperelliptic and that its canonical model is contained in a smooth quadric $Q\subset \P^3$ which is split over $k'$. This means that the character $\varepsilon^{} _Q$ is trivial on the sugroup $\mathrm{Gal}(\bar k /k')$ of $\gal$; in other words, $\varepsilon^{} _Q$ is either trivial or equal to $\varepsilon^{}_A$.
\par  It remains to prove that  $\varepsilon^{} _Q$ is nontrivial when $k'\neq k$, that is,  the nontrivial automorphism $\sigma $ of $k'/k$ exchanges the two rulings of $Q$, or equivalently the two $g^1_3$ of $X_a$. 
\par  We have ${}^\sigma \varphi =-\varphi $ (\ref{explain}).
We write as before $\varphi (\Theta)=\Sym^2C-D$; we observe that 
${}^\sigma (\varphi (\Theta ))=-\varphi (\Theta )$, hence ${}^\sigma D\sim K-D$. 
Recall that the maps $t:X_a\rightarrow |K+\varphi (a)|$ and $t':X_a\rightarrow |K-\varphi (a)|$ defining the two $g^1_3$ are given by
\begin{eqnarray*}
t(\bar z)&=& \lr{D+\varphi (z)}+\lr{K-D-\varphi (z)+\varphi (a)} \\
t'(\bar z)&=& \lr{K-D-\varphi (z)}+\lr{D+\varphi (z)-\varphi (a)}\ ,
\end{eqnarray*}where $z$ is a point of $ \tilde X_a$ and $\bar z$ its  image  in $X_a$.

Using ${}^\sigma \varphi =-\varphi $ and ${}^\sigma D\sim K-D$ we get
\[{}^\sigma t(\bar z)=\lr{K-D-\varphi (z)}+\lr{D+\varphi (z)-\varphi (a)}=t'(\bar z)\ ;\]
thus $\sigma $ exchanges $t$ and $t'$, hence the two rulings of $Q$.\qed

\bigskip

\par  One can describe  the extension $k'/k$ (hence the character $\varepsilon^{} _Q$)  using the even Clifford algebra $C^+(Q)$ \cite{De}: its center  is isomorphic to $k'$ if $k'\neq k$ and to $k\times k$ otherwise. From the description of this center  (see \cite{Bo}, \S 9, no. 4, Remarque 2), we obtain:
\begin{prop}
Assume $\mathrm{char}(k)\neq 2$, and let $\delta \in k^*$ be the discriminant of $Q$  (well defined \emph{mod.} $k^{*2}$). The extension $k'$ is isomorphic to $k(\sqrt{\delta })$.
\end{prop}
\par  Similarly, if $\mathrm{char}(k)= 2$, we have $k'=k(\lambda )$ with $\lambda ^2+\lambda =\Delta $, where $\Delta $ is the \emph{pseudo-discriminant} of $Q$ (\cite{Bo}, \S 9, exerc. 9). 
\smallskip
\par  Finally let us observe that the existence of a rational theta divisor is automatic when $k$ is finite. Indeed the theta divisors in the class $\theta $ form a torsor under $A$; by a theorem of Lang \cite{L},  such a torsor is trivial.

\end{document}